\documentclass{article}%
\usepackage{amsthm,amssymb,latexsym,amsmath}
\usepackage{amsmath}
\usepackage{amsfonts}
\usepackage{amssymb}
\usepackage{graphicx}%
\input{xypic.tex} \xyoption{all}

\newcommand{\C}{\mathbb{C}}

\newcommand{\pp}{\mathbb{P}}
\newcommand{\ce}{{\cal E}}

\newcommand{\id}{\mathbf{1}}
\newtheorem*{mthm}{Main Theorem}
\newtheorem{theorem}{Theorem}
\newtheorem{proposition}[theorem]{Proposition}

\newtheorem{remark}[theorem]{Remark}

\newtheorem*{definition}{{\bf Definition}}

\begin{document}

\title{Nahm Transform for Higgs bundles}
\author{Pedro Frejlich and Marcos Jardim\\IMECC - UNICAMP \\Departamento de Matem\'atica \\ Caixa Postal 6065 \\13083-970 Campinas-SP, Brazil}
\maketitle

\begin{abstract}
We construct the Nahm transform for Higgs bundles over a Riemann
surface of genus at least 2 as hyperholomorphic connections on the
total space of the tangent bundle of its dual Jacobian.


\end{abstract}


\section{Introduction}

Roughly speaking, the Nahm transform is a nonlinear analogue of the
Fourier transform, transforming anti-self-dual connections on the
Euclidean $\mathbb{R}^{4}$ which are invariant under a subgroup of
translations $\Lambda\subset\mathbb{R}^{4}$ into anti-self-dual
connections on the dual Euclidean space $(\mathbb{R}^{4})^{\vee}$
which are invariant under the dual subgroup of translations
$\Lambda^{\vee}\subset(\mathbb{R}^{4})^{\vee}$. Although the
construction is in principle well understood for any subgroup of
translations $\Lambda$ (see \cite{J5} for a survey), analytical
details vary for each $\Lambda$.

Particularly relevant to this paper is the case $\Lambda=\mathbb{R}^{2}%
\times\mathbb{Z}^{2}$, $\Lambda^{\vee}=\mathbb{Z}^{2}$, considered by the
second named author in \cite{J1,J2} and in collaboration with O. Biquard in \cite{BJ}.
There, it was shown that the Nahm transform provides a 1-1 correspondence between singular solutions of
Hitchin's equations on a 2-dimensional torus $T^{2}$ and anti-self-dual
connections on $(T^{2})^{\vee}\times\mathbb{R}^{2}$ satisfying certain
conditions at infinity (which are equivalent to the existence of a holomorphic
extension to $(T^{2})^{\vee}\times S^{2}$). Here, $(T^{2})^{\vee}$ is meant to
denote the 2-dimensional torus dual to $T^{2}$. The main goal of the present
paper is to generalize the Nahm transform of solutions of Hitchin's equations
on a 2-dimensional torus $T^{2}$ to solutions of Hitchin's equations on a
Riemann surface $\Sigma$ of genus $g\ge2$. More precisely, denoting by $J^{\vee}$ 
the dual to the Jacobian of $\Sigma$ we prove:

\begin{mthm}
The Nahm transform of an irreducible solution of Hitchin's equations of rank at least 2
over a Riemann surface $\Sigma$ of genus $g\ge2$ is a Hermitian vector bundle
$\widehat{\mathcal{E}}\to J^{\vee}\times H^{0}(K_{\Sigma})$ equipped with a
hyperholomorphic unitary connection $\widehat\nabla$.
Moreover, the holomorphic structure induced by $\widehat\nabla$ on $\widehat{\mathcal{E}}$
with respect to a product complex structure on $J^{\vee}\times H^{0}(K_{\Sigma})$ extends to a
holomorphic bundle over $J^{\vee}\times\mathbb{P}(H^{0}(K_{\Sigma})\oplus\mathbb{C})$.
\end{mthm}

The problem of characterizing the Nahm transform of solutions of Hitchin's
equations has a long history. It was first suggested in
Garc\'{\i}a-Prada's PhD thesis \cite[Section 4.4]{GP} that the Nahm transform
of a unitary Hermitian-Einstein connection on $\Sigma$ should be a
Hermitian-Einstein connection on its Jacobian. More recently, Tejero Prieto
considered the Nahm and Fourier--Mukai transforms for Riemann surfaces of higher
genus \cite{TP}. Also closely related to the present paper is Bonsdorff's Fourier--Mukai
transform for stable Higgs bundles on $\Sigma$ \cite{Bo}; in fact our construction,
which is carried out in Section \ref{nt}, is the differential geometric analogue of the algebraic
geometric construction of \cite{Bo}, with the advantage of providing one additional
piece of information: the existence of the hyperholomorphic unitary connection on
$J^{\vee}\times H^{0}(K_{\Sigma})$; using twistor theory, Bonsdorff constructed an auto-dual
connection, see \cite{Bo-t}.

This paper also completes the picture for the Nahm transform of Higgs bundles on Riemann surfaces.
The genus $0$ case was recently considered by Szabo in \cite{S}, while the genus $1$ case was studied
in \cite{J1,J2}, as mentioned above. One peculiar aspect of the story is that, since there are
no smooth Higgs bundles over surfaces of genus $0$ and $1$, the authors of \cite{J1,J2,S} were forced to
consider {\em singular} Higgs bundles. In this work we only consider the Nahm transform of smooth Higgs bundles;
an interesting sequel to this paper would be the construction of the Nahm transform for singular Higgs bundles over
surfaces of genus at least two, and see what new features arise. In addition, analytical aspects of the transformed hyperholomorphic connection which were studied in the genus $0$ and genus $1$ cases in \cite{BJ,J1,J2,S} are not explored in this paper and deserve further research.

\bigskip

\paragraph{Acknowledgements.}
P.F.'s research was supported by the CNPq grant number 130226/2005-0. M.J. is
partially supported by the CNPq grant number 300991/2004-5, the FAPESP grant 2005/04558-0,
and FAEPEX. We thank Benoit Charbonneau for his comments on a preliminary version of this paper.


\section{Riemann surfaces, Jacobians and the Poincar\'e line bundle}\label{s}

Let $\Sigma$ be a closed Riemann surface of genus $g>1$, regarded as the
closed, oriented (real) surface obtained from a regular $4g$-gon $\Delta$
(whose edges we label as $\alpha_{1},\beta_{1},\alpha_{1}^{-1},\beta_{1}%
^{-1},...,\beta_{g}^{-1}$ with gluing relation $\Pi_{i}[\alpha_{i},\beta
_{i}]=1$) endowed with a conformal class of Riemannian metrics $\mathcal{G}$.
We denote by $a_{i}$ and by $b_{i}$ the loops determined by the $\alpha_{i}$'s
and $\beta_{i}$'s, respectively. The integral homology of $\Sigma$ in
dimension $1$ is then spanned by the homology classes determined by $a_{i}$ and $b_{i}$.

For dimensional reasons, every metric on such $\Sigma$ has closed
$(1,1)$-form and is therefore K\"{a}hler. We thus obtain decompositions (cf.
\cite[Ch. 0]{GH})
\[
H^{n}\left(  \Sigma;\mathbb{C}\right)  =\bigoplus\limits_{p+q=n}H^{p,q}\left(
\Sigma;\mathbb{C}\right)
\]%
\[
H^{p,q}\left(  \Sigma;\mathbb{C}\right)  =\overline{H^{q,p}\left(
\Sigma;\mathbb{C}\right)  }%
\]
(where $n=0,1,2$, $p,q=0,1$ and $H^{p,q}\left(  \Sigma;\mathbb{C}\right)  $
stands for the Dolbeault cohomology) which hold throughout the conformal class
$\mathcal{G}$. In particular
\[
H^{1}\left(  \Sigma;\mathbb{C}\right)  =H^{1,0}\left(  \Sigma;\mathbb{C}%
\right)  \oplus\overline{H^{1,0}\left(  \Sigma;\mathbb{C}\right)  }%
\]
so that, fixing a basis $\omega_{1},...,\omega_{g}$ for $H^{0}(K_{\Sigma})$,
de Rham's theorem shows that if we set $A_{i}\left(  \omega_{j}\right)  =\int_{a_{i}}\omega_{j}$ and $B_{i}\left(  \omega_{j}\right)  =\int_{b_{i}}\omega_{j}$, the $2g$ elements
\[
A_{i}=%
\begin{pmatrix}
A_{i}(\omega_{1})\\
\vdots\\
A_{i}(\omega_{g})
\end{pmatrix}
~~,~~B_{i}=%
\begin{pmatrix}
B_{i}(\omega_{1})\\
\vdots\\
B_{i}(\omega_{g})
\end{pmatrix}
\]
determine an $\mathbb{R}$-basis $\Lambda$ for the vector space
$\mathbb{C}^{g}$. Thus $J=\mathbb{C}^{g}/\Lambda$, called the \emph{Jacobian} of the Riemann surface
$\Sigma$, is a complex torus of dimension $g$ which clearly depends only on the conformal structure
$\mathcal{G}$ on $\Sigma$.

Once we fix a base point $p\in\Sigma$, we may define Abel's map $\mathsf{Ab}%
_{p}:\Sigma\longrightarrow J$ by:
\[
\mathsf{Ab}_{p}(q)=%
\begin{pmatrix}
\int_{p}^{q}\omega_{1}\\
\vdots\\
\int_{p}^{q}\omega_{g}%
\end{pmatrix}
\]
where $\int_{p}^{q}$ denotes integration along any piecewise smooth
path joining $p$ to $q$. This (holomorphic)\ mapping can be naturally extended
to a mapping of $\mathsf{Div}^{0}(\Sigma)$ (the group of degree zero divisors
on $\Sigma$) onto $J$:
\[
\mathsf{Ab}_{p}:\mathsf{Div}^{0}(\Sigma)\longrightarrow J~~.
\]
By means of Riemann's bilinearity relations (\cite[Ch. 2]{GH}) one can show that any
such divisor which is mapped to $0$ by $\mathsf{Ab}_{p}$ is associated to a
global meromorphic function, and thus $J$ can be regarded as the quotient of
$\mathsf{Div}^{0}(\Sigma)$ modulo principal, degree-zero divisors:
\[
J=\frac{\mathsf{Div}^{0}(\Sigma)}{\mathsf{PDiv}^{0}(\Sigma)}~~.
\]
This leads to an identification between $J$ and $\mathsf{Pic}^{0}(\Sigma)$, the
\emph{Picard group} of topologically trivial, holomorphic line bundles on
$\Sigma$. Alternatively, a point in $J$ can also be thought of as a flat
connection on the topologically trivial line bundle $\underline{\mathbb{C}}\rightarrow\Sigma$. 

Once and for all, we fix (a) the standard flat metric on $J$, (b) a base point $p\in\Sigma$
and (c) the corresponding metric on $\Sigma$ induced by the Abel map $\mathsf{Ab}_{p}$.

\begin{proposition}
\label{spin} Given a spin structure on $J$, there exists a spin structure on
$\Sigma$ which is compatible with the Abel map $\mathsf{Ab}_{p}:\Sigma
\hookrightarrow J$.
\end{proposition}

\begin{proof}
The key point in this proposition is that a spin submanifold of a manifold endowed with a spin structure has a canonical spin structure
\textit{once we choose a spin structure for its normal bundle} \cite[Chapter 2, Proposition 1.2]{LM}.
Since the normal bundle $N\Sigma\longrightarrow\Sigma$ is oriented and spin, it is (differentially) trivial whenever its (real) rank is at least three, which happens precisely when $g>2$, see \cite[Chapter 2, Proposition 2 2.15]{LM}. In that case, it is obvious that we can endow $N\Sigma$ with a spin structure, and therefore $T\Sigma$ is assigned a spin structure compatible with that of $J$.
As for the case $g=2$, we remark that $c_{1}\left(  J\right)  =0$ whence
\[
0=\left(  \mathsf{Ab}_{p}\right)  ^{\ast}c_{1}\left(  J\right)  =c_{1}\left(
N\Sigma\oplus T\Sigma\right)  =c_{1}\left(  N\Sigma\right)  +c_{1}\left(
T\Sigma\right)
\]
and therefore $N\Sigma$, regarded as an $U_1$-bundle, is differentially isomorphic to the inverse of the tangent bundle of $\Sigma$.
Next we recall that a spin structure on a $U_1$-bundle $P$ on $\Sigma$ is simply a bundle map
$\widetilde{P}\longrightarrow P$ which restricts to the two-sheeted covering $z\longmapsto z^{2}$ on each fiber. So, if $P_{U_{1}}\left(  T\Sigma\right)  $ is given by the \v{C}ech cocycle
\[
g_{\alpha\beta}:U_{\alpha}\cap U_{\beta}\longrightarrow U_{1}%
\]
a spin structure
\[
\widetilde{P}\longrightarrow P_{U_{1}}\left(  T\Sigma\right)
\]
is represented by $\left\{  h_{\alpha\beta}\right\}  $ with $h_{\alpha
\beta}^{2}=g_{\alpha\beta}$; thus $P_{U_{1}}\left(  T\Sigma^{-1}\right)  $ is represented by
$\{  g_{\alpha\beta}^{-1}\}$ so that $\{  h_{\alpha\beta}%
^{-1}\}  $ is a \v{C}ech cocycle which determines a spin structure on $T\Sigma^{-1}\simeq N\Sigma$ :
\[
\widetilde{Q}\longrightarrow P_{U_{1}}\left(  \left(  T\Sigma\right)
^{-1}\right)
\]
Therefore, $N\Sigma$ can always be given a spin structure.
\end{proof}

\begin{remark}\rm
In view of this Proposition, and bearing in mind that $\Sigma$ is by construction a Riemannian submanifold of $J$,
we see that the Abel map $\mathsf{Ab}_{p}$ intertwines Clifford multiplications $\cdot_{J}$ on $J$ and $\cdot_{\Sigma}$ on $\Sigma$. By this we mean that given differential forms $\varphi,\psi$ on $J$, the following relation holds  :%
$$ \mathsf{Ab}_{p}^{\ast}\left(  \varphi\cdot_{J}\psi\right)  =\mathsf{Ab}%
_{p}^{\ast}\left(  \varphi\right)  \cdot_{\Sigma}\mathsf{Ab}_{p}^{\ast}\left(
\psi\right) $$
\end{remark}

A \emph{dual Jacobian} $J^{\vee}$ of $\Sigma$ may be described in the
following terms: to the maximal lattice $\Lambda$ we associate the dual
lattice
\[
\Lambda^{\vee}=\mathsf{Hom}\left(  \Lambda,\mathbb{Z}\right)  \subset\left(
\mathbb{C}^{g}\right)  ^{\vee}
\]
consisting of all linear functionals on $\mathbb{C}^{g}$ that assume only
integer values on $\Lambda$. The lattice $\Lambda^{\vee}$ is easily seen to be
maximal, so that the quotient $J^{\vee}=\left(\mathbb{C}^{g}\right)^{\vee}/\Lambda^{\vee}$ is also a complex torus of dimension $g$.

Since $J$ has trivial tangent bundle, any functional $\eta\in\left(\mathbb{C}^{g}\right)^{\vee}$ may be regarded as a flat connection $1$-form on the trivial line bundle $\underline{\mathbb{C}}\to J$; the induced topologically trivial, holomorphic line bundle on $J$ is denoted by $L_{\eta}$. Moreover, any such bundle is of the form $L_{\eta}$ for some $\eta$, and $L_{\eta}$ is isomorphic to $L_{\eta^{\prime}}$ if and only if $\eta$ and $\eta^{\prime}$ differ by an element of the dual lattice $\Lambda^{\vee}$; cf. \cite[Ch. 3]{DK}. Dually, any topologically trivial, holomorphic line bundle $L\longrightarrow J^{\vee}$ is of the form $L_{\xi}$ for some $\xi\in\mathbb{C}^{g}$, and two such bundles are isomorphic precisely when the functionals from which they are constructed have the same image in $J$.

These remarks suggest the following construction: on the product $J\times\left(\mathbb{C}^{g}\right)^{\vee}$,
we let $\Lambda^{\vee}$ act on the trivial line bundle
$\underline{\mathbb{C}}=J\times\left(\mathbb{C}^{g}\right)^{\vee}\times\mathbb{C}\rightarrow
J\times\left(\mathbb{C}^{g}\right)^{\vee}$ with product connection by the formula
\begin{align*}
\Lambda^{\vee}\times\underline{\mathbb{C}}  &  \longrightarrow\underline
{\mathbb{C}}\\
\left(  \mu,\left(  \xi,\eta,z\right)  \right)   &  \longmapsto\left(
\xi,\eta+\mu,e^{2\pi i\mu(\xi)}z\right)
\end{align*}
The quotient of $\underline{\mathbb{C}}$ by this action preserves the connection and defines a line bundle $\mathcal{P}\longrightarrow J\times J^{\vee}$ with connection $1$-form $\omega$ which, in terms of dual, orthonormal real local coordinates $\{\xi_{j}\}_{j=1}^{2g}$ of $J$ and $\{\eta_{j}\}_{j=1}^{2g}$ of $J^{\vee}$ we can write as
\begin{equation}
\omega=2\pi i\sum_{j=1}^{2g}\left(  \eta_{j}d\xi_{j}-\xi_{j}d\eta_{j}\right)
\label{pconn}%
\end{equation}
This is called the \emph{Poincar\'{e} line bundle} of $J\times J^{\vee}$; accordingly, the connection with which it is endowed is called the Poincar\'{e} connection. Its remarkable feature is that it provides a \emph{duality data} for topologically trivial, holomorphic line bundles
between $J$ and $J^{\vee}$, i.e.:
\begin{align*}
\mathcal{P}|_{J\times\{\eta\}}  &  \simeq L_{\eta}\longrightarrow J\\
\mathcal{P}|_{\{\xi\}\times J^{\vee}}  &  \simeq L_{-\xi}\longrightarrow
J^{\vee}%
\end{align*}
Notice that the curvature of the Poincar\'{e} connection (\ref{pconn}) is
given by
\begin{equation}
\Omega_{\mathcal{P}}=d\omega=4\pi i\sum_{j=1}^{2g}d\eta_{j}\wedge d\xi_{j}~~; \label{pcurv}%
\end{equation}
Choosing dual unitary bases
$\left\{\gamma_{k}=\xi_{2k-1}+i\xi_{2k}\right\}_{k=1}^{g}$ and $\left\{\rho_{k}=\eta_{2k-1}+i\eta_{2k}\right\}_{k=1}^{g}$ for $\mathbb{C}^{g}$ and $(\mathbb{C}^{g})^{\vee}$, respectively, we get:
$$ \Omega_{\mathcal{P}} = 
2\pi i \sum_{k=1}^{g}(d\rho_{k}\wedge d\overline{\gamma}_{k}+d\overline{\rho}_{k}\wedge d\gamma_{k}) $$
In particular, $\Omega_{\mathcal{P}}$ is of type $(1,1)$ as a 2-form on $J\times J^{\vee}$ with respect to the product complex structure, and as such it induces a holomorphic structure on $\mathcal{P}$.

In what follows, $K_{\Sigma}$ will (as usual) denote the canonical bundle of $\Sigma$.
Recall that the space of abelian differentials $H^{0}(K_{\Sigma})$ has dimension $g$; this allows us to provide 
$J^{\vee}\times H^{0}(K_{\Sigma})$ with a flat hyperk\"ahler structure, as follows. As above, let $\{\eta_{j}\}_{j=1}^{2g}$ be orthonormal real local coordinates of $J^\vee$ and let $\{s_{j}\}_{j=1}^{2g}$ be a real basis of $H^{0}(K_{\Sigma})$ so that $\{\sigma_k=s_{2k-1}+is_{2k}\}_{k=1}^{g}$ is a basis of holomorphic 1-forms on $\Sigma$. We define the following complex structures on the product $J^{\vee}\times H^{0}(K_{\Sigma})$
$$ \begin{array}{lcr}
I_1(\partial/\partial\eta_{2k-1})=\partial/\partial\eta_{2k} &
I_2(\partial/\partial\eta_{2k-1})=\partial/\partial s_{2k-1} &
I_3(\partial/\partial\eta_{2k-1})=\partial/\partial s_{2k} \\
I_1(\partial/\partial\eta_{2k})=-\partial/\partial\eta_{2k-1} &
I_2(\partial/\partial\eta_{2k})=-\partial/\partial s_{2k} &
I_3(\partial/\partial\eta_{2k})=\partial/\partial s_{2k-1} \\
I_1(\partial/\partial s_{2k-1})=\partial/\partial s_{2k} &
I_2(\partial/\partial s_{2k-1})=-\partial/\partial\eta_{2k-1} &
I_3(\partial/\partial s_{2k-1})=-\partial/\partial\eta_{2k} \\
I_1(\partial/\partial s_{2k})=-\partial/\partial s_{2k-1} &
I_2(\partial/\partial s_{2k})=\partial/\partial\eta_{2k} &
I_3(\partial/\partial s_{2k})=-\partial/\partial\eta_{2k-1}
\end{array} $$
The flatness of the metric on $J^{\vee}\times H^{0}(K_{\Sigma})$ makes these structures integrable. Note that $I_1$ is the product of a complex structure on $J^{\vee}$ with a complex structure on $H^{0}(K_{\Sigma})$.

Finally, notice also that $J^{\vee}\times H^{0}(K_{\Sigma})$ can be identified with the total space of the tangent bundle of $J^{\vee}$.


\section{Higgs bundles and hyperholomorphic connections} \label{e}

Let $\mathcal{E}\to\Sigma$ be a Hermitian vector bundle provided with a unitary
connection $\nabla$ and a bundle morphism $\theta:\mathcal{E}\to \mathcal{E}\otimes K_{\Sigma}%
$. Hitchin's equations for $\nabla$ and $\theta$
\begin{equation}
\label{hit}\left\{
\begin{array}
[c]{c}%
F_{\nabla}+ [\theta,\theta^{*}] = 0\\
\overline{\partial}_{\nabla}\theta= 0
\end{array}
\right.
\end{equation}
were introduced in \cite{H} as the dimensional reduction of the
anti-self-duality equations from 4 to 2 dimensions, and have been studied
extensively in the past two decades. Relevant to us is the fact that $\nabla$
induces a holomorphic structure on the bundle $\mathcal{E}$, so that the second equation
simply says that the morphism $\theta$ is holomorphic with respect to the connection.

The pair $(\mathcal{E},\theta)$ consisting of a holomorphic bundle $\mathcal{E}$ and a section
$\theta\in H^{0}(\mathrm{End}(\mathcal{E})\otimes K_{\Sigma})$ is known as a \emph{Higgs bundle};
$\theta$ is called the \emph{Higgs field}. The trivial Higgs bundle is the pair $(\mathcal{O}_{\Sigma},0)$.
A morphism between Higgs bundles $(\mathcal{E},\theta)$ and $(\mathcal{E}^{\prime},\theta^{\prime})$ is simply
a bundle morphism $\phi:\mathcal{E}\to\mathcal{E}^{\prime}$ rendering the following diagram commutative:
\[
\xymatrix{
\ce \ar[d]_{\phi} \ar[r]^{\theta} & \ce\otimes K_\Sigma \ar[d]^{\phi \otimes 1} \\
\ce^{\prime} \ar[r]^{\theta^{\prime}} & \ce^{\prime}\otimes K_\Sigma
}
\]
Accordingly, a Higgs subbundle of $(\mathcal{E},\theta)$ is merely a subbundle
$\mathcal{E}^{\prime}$ of $\mathcal{E}$ which $\theta$ maps into
$\mathcal{E}^{\prime}\otimes K_{\Sigma}$.

The bundle $(\mathcal{E},\theta)$ is said to be \emph{stable} if the \emph{slope} $\mu(\mathcal{E})=\mathrm{deg}(\mathcal{E}%
)/\mathrm{rk}(\mathcal{E})$ of $\mathcal{E}$ is strictly greater
than that of any of its proper Higgs subbundles;
$(\mathcal{E},\theta)$ is called \emph{semistable} if equality can
occur, and \emph{polystable} if it is the direct sum of stable Higgs
bundles with the same slope.

We also may regard a Higgs bundle $(\mathcal{E},\theta)$ as a 2-term complex
of coherent $\mathcal{O}_{\Sigma}$-modules:
\[
\mathbf{E} = \left\{  \mathcal{E}\overset{\theta}{\longrightarrow}
\mathcal{E}\otimes K_{\Sigma}\right\}  ~~.
\]
This allows us to compute the hypercohomology $\mathbb{H}^{\ast}\left(\mathbf{E}\right)$, as described in \cite[Ch. 3]{GH}.

A key ingredient to the Nahm transform we are about to define will be the following vanishing theorem due to Hausel \cite[Corollary 5.1.4]{Ha}:

\begin{proposition}\label{v1}
If $\mathbf{E}=(\mathcal{E},\theta)$ is a nontrivial stable Higgs bundle with zero slope,
then both $\mathbb{H}^{0}\left(\mathbf{E}\right)$ and $\mathbb{H}^{2}\left(\mathbf{E}\right)$ vanish.
\end{proposition}

This is essentially a consequence of (1) a famous lemma by Narasimhan and Seshadri \cite{NS}, which states that if $f:E\longrightarrow F$ is a nontrivial morphism between holomorphic vector bundles over a Riemann surface $\Sigma$, then there exist holomorphic vector bundles $E_{1},F_{1},E_{2}$ and $F_{2}$ and a morphism $E_{2}\longrightarrow F_{1}$ which is an isomorphism over an open Zariski subset $U$ of $\Sigma$, fitting into the commutative diagram with exact rows

\[
\xymatrix{0 \ar[r] & E_1 \ar[r] & E \ar[d]^f \ar[r] & E_2 \ar[d] \ar[r] & 0 \\
0 & F_2 \ar[l] & F \ar[l] & F_1 \ar[l] & 0 \ar[l]
}
\]
and (2) the fact that $E_{2}\longrightarrow F_{1}$ being an isomorphism over $U$ implies $\mathrm{deg}(F_{1})\ge\mathrm{deg}(E_{2})$.

We are now in position to state the fundamental result first proved by Hitchin in \cite{H} and later generalized by other authors:

\begin{theorem}
\label{hit.thm} If $(\nabla,\theta)$ is an irreducible solution of Hitchin's
equations (\ref{hit}) on a complex vector bundle $E$, then the associated
Higgs bundle is stable of degree zero. Conversely, to each stable Higgs bundle
$(\mathcal{E},\theta)$ of degree zero, there corresponds an irreducible
unitary connection $\nabla$, unique up to gauge equivalence, which is
compatible with the holomorphic structure of $\mathcal{E}$ and satisfies the
first equation of (\ref{hit}).
\end{theorem}

Finally, let $(X,\{I,J,K\})$ be a hyperk\"ahler manifold. Given real numbers $a,b$ and $c$ with $a^{2}+b^{2}+c^{2}=1$, the complex structure $L=aI+bJ+cK$ is said to be \emph{induced} by the hyperk\"ahler structure on $X$. A connection
$\nabla$ on a Hermitian vector bundle $V$ over $X$ is said to be \emph{hyperholomorphic} if its curvature $F_{\nabla}$ is of type $(1,1)$ with respect to \emph{any} complex structure $L$ induced by the hyperk\"ahler structure.
The Poincar\'e connection (\ref{pconn}) is an example of a hyperholomorphic connection, as is shown in
\cite[Theorem 4]{BaJ}. It also is interesting to note that if $\dim_{\mathbb{R}}X=4$, then a connection is anti-self-dual if and only if it is hyperholomorphic, so that hyperholomorphicity may be thought of as a generalization of the anti-self-dual condition in four dimensions.

If we think of $SU(2)$ as the space of unit quaternions, we obtain an obvious action of $SU(2)$ on the bundle of differential forms over $X$, which one can easily check to be parallel, and therefore commutes with the trace Laplacian $\Delta$ on $X$. By Hodge theory \cite[Ch. 0]{GH} it follows that if $X$ is compact, this natural action of $SU(2)$ descends to an action on the cohomology level. This is the point where two results of Verbitsky \cite{V} come into play : the first states that $\omega\in H^{2p}(X;\mathbb{C})$ is fixed by this action if and only if $\omega\in H_{L}^{p,p}(X)$ for every induced complex structure $L$, whereas the second asserts that the subspace $H^{inv}\subset H^{2}(X;\mathbb{C})$ of $SU(2)$-invariant $2$-forms is orthogonal to the $3$-dimensional subspace generated by the K\"{a}hler forms $\omega_{I},\omega_{J},\omega_{K}$ (\cite[Lemma 2.1]{V}). Consequently, every hyperholomorphic connection $\nabla$ on $V$ is Hermitian-Einstein with constant zero with respect to any complex structure $L$ induced by the hyperk\"ahler structure; in other words, a hyperholomorphic connection $\nabla$ satisfies $\Lambda_{L}F_{\nabla}= 0$, where $\Lambda_{L}$ denotes contraction with the K\"ahler form associated with the complex structure $L$.


\section{The transform} \label{nt}

Let $\mathbf{E}=(\mathcal{E},\theta)$ be a stable Higgs bundle over $\Sigma$
of rank at least two, and let $\nabla$ be the associated unitary connection,
in the sense of Theorem \ref{hit.thm}. According to our previous remarks on
the Poincar\'{e} bundle, given $\eta\in J^{\vee}$ we have a topologically
trivial, holomorphic line bundle $L_{\eta}|_{\Sigma}\to\Sigma$.

Now to each pair $(\eta,\sigma)\in J^{\vee}\times H^{0}(K_{\Sigma})$ there
corresponds a Higgs bundle $\mathbf{E}(\eta,\sigma) = (\mathcal{E}(\eta),\theta_{(\eta,\sigma)})$, where
$$ \mathcal{E}(\eta) = \mathcal{E}\otimes\left(  L_{\eta}|_{\Sigma}\right) ~\mathrm{and}~
\theta_{(\eta,\sigma)} = \theta\otimes\mathbf{1}_{L_{\eta}}+\mathbf{1}_{\mathcal{E}(\eta)}\otimes\sigma~~, $$
where $\sigma\in H^{0}(K_{\Sigma})$ is now regarded as a map $\sigma:\mathcal{O}_{\Sigma}\to K_{\Sigma}$.
Clearly, every such Higgs bundle is stable, since $(\mathcal{E},\theta)$ is. Furthermore, the degree of each
$\mathcal{E}(\eta)$ is zero, since $L_{\eta}$ is topologically trivial.

The unitary connection associated to $(\mathcal{E}(\eta),\theta_{(\eta,\sigma)})$
is the tensor connection $\nabla_{\eta}=\nabla\otimes\mathbf{1}_{L_{\eta}}+\mathbf{1}_{E}\otimes2\pi i\eta$,
and it is not difficult to see that, for each $(\eta,\sigma)\in J^{\vee}\times H^{0}(K_{\Sigma})$, the
pair $(\nabla_{\eta},\theta_{(\eta,\sigma)})$ satisfies Hitchin's equations
(\ref{hit}). Notice that $\overline{\partial_{\nabla_{\eta}}}$ yields the
holomorphic structure of $\mathcal{E}(\eta)$. To simplify notation, we use
$\overline{\partial}_{\eta}$ to denote $\overline{\partial_{\nabla_{\eta}}}$.

With this in mind, we define a family of first order differential operators
parameterized by $J^{\vee}\times H^{0}(K_{\Sigma})$ in the following manner:
$$ \mathcal{D}_{(\eta,\sigma)}~:~
L_{p+1}^{2}\left(\mathcal{E}(\eta)\otimes(\Lambda_{\Sigma}^{0}\oplus\Lambda_{\Sigma}^{1,1})\right)
\longrightarrow
L_{p}^{2}\left(\mathcal{E}(\eta)\otimes(\Lambda_{\Sigma}^{1,0}\oplus\Lambda_{\Sigma}^{0,1})\right) $$
$$ \mathcal{D}_{(\eta,\sigma)} = -\left(
\overline{\partial}_{\eta}+\overline{\partial}_{\eta}^{\ast}\right) +
\left(  \theta_{(\eta,\sigma)}-\theta_{(\eta,\sigma)}^{\vee}\right) $$
where $\theta^{\vee}=i\theta^*\circ\Lambda+i\Lambda\circ\theta^*$ with $\Lambda$ denoting the
contraction with the K\"ahler form on $\Sigma$. More precisely, setting
$\psi_{0}\in L^{2}(\mathcal{E}(\eta)\otimes\Lambda_{\Sigma}^{0})$ and
$\psi_{2}\in L^{2}(\mathcal{E}(\eta)\otimes\Lambda_{\Sigma}^{1,1})$, one has:
$$ \mathcal{D}_{(\eta,\sigma)}(\psi_{0},\psi_{2})=
(-\overline{\partial}_{\eta}^{\ast}\psi_{2}+\theta_{(\eta,\sigma)}\psi_{0},
-\overline{\partial}_{\eta}\psi_{0} -\theta_{(\eta,\sigma)}^{\vee}\psi_{2})~~. $$
Modulo the algebraic term $\left(  \theta_{(\eta,\sigma)}-\theta_{(\eta,\sigma)}^{\vee}\right)$,
the operator $\mathcal{D}_{(\eta,\sigma)}$ is just the canonical Dirac operator of a Hermitian manifold;
$\mathcal{D}_{(\eta,\sigma)}$ is therefore an elliptic differential operator for each
choice of $(\eta,\sigma)$.

\begin{remark}\rm
In \cite{H2}, Hitchin introduced the following Dirac operator for Higgs bundles:
$$ D ~:~ L_{p+1}^{2}\left(\mathcal{E}\otimes(\Lambda_{\Sigma}^{0}\oplus\Lambda_{\Sigma}^{0})\right)
\longrightarrow
L_{p}^{2}\left(\mathcal{E}\otimes(\Lambda_{\Sigma}^{1,0}\oplus\Lambda_{\Sigma}^{0,1})\right) $$
$$ D(\phi_1,\phi_2) = (-\partial_\nabla\phi_1+\theta\phi_2,-\theta^*\phi_1-\overline{\partial}_\nabla\phi_2) ~~. $$
Using the K\"ahler identity $\overline{\partial}_{\nabla}^{\ast}=i{\partial}_{\nabla}\Lambda$ on $(1,1)$-forms,
it is easy to see that:
$$ {\cal D}_{(0,0)}(\psi_0,\psi_2) = D(i\Lambda\psi_2,\psi_0) ~~.$$
\end{remark}

\begin{proposition}
\label{ker=0} If $(\mathcal{E},\theta)$ is a stable Higgs bundle of degree
zero and rank at least two, then $\ker\mathcal{D}_{(\eta,\sigma)}=0$ for each
$(\eta,\sigma)$.
\end{proposition}

\begin{proof}
First notice that $\mathcal{D}_{(\eta,\sigma)}$ unfolds as a three-term complex of vector spaces
\begin{equation}\label{c1} \mathbf{M}_{(\eta,\sigma)} = \{ \xymatrix{
L^2_{p+1}(\mathcal{E}(\eta)) \ar[r]^{\!\!\!\!\!\!\!\!\!\!\!\!\!\!\!\!\!\!\!\!-\overline{\partial}_{\eta}+\theta_{(\eta,\sigma)}} &
L^2_p(\mathcal{E}(\eta)\otimes(\Lambda_{\Sigma}^{1,0}\oplus\Lambda_{\Sigma}^{0,1}))
\ar[r]^{~~-\overline{\partial}_{\eta}-\theta_{(\eta,\sigma)}} & L^2_{p-1}(\mathcal{E}(\eta)\otimes\Lambda^{1,1}_{\Sigma})
} \} \end{equation}
so that
$$ \ker\mathcal{D}_{(\eta,\sigma)} = H^0(\mathbf{M_{(\eta,\sigma)}}) \oplus H^2(\mathbf{M_{(\eta,\sigma)}}) ~~{\rm and} $$
$$ {\rm coker}~\mathcal{D}_{(\eta,\sigma)} = H^1(\mathbf{M_{(\eta,\sigma)}}) ~~. $$
Using Hodge theory, it can be shown that $H^q(\mathbf{M_{(\eta,\sigma)}})$ is canonically isomorphic to the
hypercohomology $\mathbb{H}^{q}$ of the twisted complex
$$ \mathbf{E}(\eta,\sigma) = \left\{  \ce(\eta)\stackrel{\theta_{(\eta,\sigma)}}{\longrightarrow}
\ce(\eta)\otimes K_{\Sigma} \right\} ~~, $$
see \cite[Section 7]{H2}. Since for each $(\eta,\sigma)\in J^{\vee}\times H^{0}(K_{\Sigma})$ the Higgs bundle
$\mathbf{E}(\eta,\sigma)$ is stable of slope zero, it follows from Proposition \ref{v1} that
$\mathbb{H}^{0}(\mathbf{E}(\eta,\sigma))$ and $\mathbb{H}^{2}(\mathbf{E}(\eta,\sigma))$
vanish, hence $\ker\mathcal{D}_{(\eta,\sigma)}=0$.
\end{proof}

We point out that we may fix a topological isomorphism
\[
\mathcal{E}\otimes\underline{\mathbb{C}}\overset{\simeq}{\longrightarrow}\mathcal{E}
\]
to identify each $L^{2}_{p}\left(\mathcal{E}(\eta)\otimes(\Lambda_{\Sigma}^{1,0}%
\oplus\Lambda_{\Sigma}^{0,1})\right)  $ with $L^{2}_{p}\left(\mathcal{E}\otimes
(\Lambda_{\Sigma}^{1,0}\oplus\Lambda_{\Sigma}^{0,1})\right)  $ and think of
$\mathrm{coker}~\mathcal{D}_{(\eta,\sigma)}$ as a subspace of $L^{2}%
_{p}\left(\mathcal{E}\otimes(\Lambda_{\Sigma}^{1,0}\oplus\Lambda_{\Sigma}%
^{0,1})\right)  $. Direct calculation shows that the dimension of $\mathrm{coker}~\mathcal{D}_{(\eta,\sigma)}$ does not vary in the $\sigma$-direction; i.e., that $\mathrm{dim}~\mathrm{coker}~\mathcal{D}_{(\eta,\sigma)}=\mathrm{dim}~\mathrm{coker}~\mathcal{D}_{(\eta,0)}$, whence
\[
J^{\vee}\times H^{0}(K_{\Sigma})\longrightarrow\mathbb{Z}
\]
\[
(\eta,\sigma)\longmapsto\mathrm{dim}~\mathrm{coker}~\mathcal{D}_{(\eta,\sigma)}
\]
factors through the projection $J^{\vee}\times H^{0}(K_{\Sigma})\longrightarrow J^{\vee}$; compactness of $J^{\vee}$ allows us then to use the usual index-theoretic machinery (see \cite[Chapter 3]{DK}) to conclude that the
dimension of $\mathrm{coker}~\mathcal{D}_{(\eta,\sigma)}$ is constant along
$J^{\vee}\times H^{0}(K_{\Sigma})$, since $\ker\mathcal{D}_{(\eta,\sigma)}$
vanishes for each $(\eta,\sigma)$. Therefore the assignment
\[
J^{\vee}\times H^{0}(K_{\Sigma}) \ni(\eta,\sigma) \longmapsto\widehat
{\mathcal{E}}_{(\eta,\sigma)}=\mathrm{coker}~\mathcal{D}_{(\eta,\sigma)}
\]
defines a vector subbundle $\widehat{\mathcal{E}}$ of the trivial Hilbert bundle
\[
H^{-} = L^{2}_{p}\left(\mathcal{E}\otimes(\Lambda_{\Sigma}^{1,0}\oplus\Lambda_{\Sigma
}^{0,1})\right)  \times J^{\vee}\times H^{0}(K_{\Sigma}) \longrightarrow
J^{\vee}\times H^{0}(K_{\Sigma}) ~~.
\]
whose direct sum with the trivial Hilbert bundle
\[
L_{p}^{2}(\mathcal{E}\otimes(\Lambda_{\Sigma}^{0}\oplus\Lambda_{\Sigma}^{1,1}))\times
J^{\vee}\times H^{0}(K_{\Sigma})\longrightarrow J^{\vee}\times H^{0}%
(K_{\Sigma})
\]
we will denote by $H$.
The rank of $\widehat{\mathcal{E}}$ is given by minus the index of the operator
$\mathcal{D}_{(0,0)}$, which is precisely the Dirac operator
\[
\overline{\partial}_{\eta}\oplus\overline{\partial}_{\eta}^{\ast} ~:~
L^{2}_{p+1}\left(\mathcal{E}(\eta)\otimes(\Lambda_{\Sigma}^{0}\oplus\Lambda_{\Sigma}^{1,1})\right)
\longrightarrow
L^{2}_{p}\left(\mathcal{E}(\eta)\otimes(\Lambda_{\Sigma}^{1,0}\oplus\Lambda_{\Sigma}^{0,1})\right)  ~~.
\]
Thus
\[
\mathrm{rk}(\widehat{\mathcal{E}}) = -2\int_{\Sigma}\mathrm{ch}(\mathcal{E})\cdot\mathrm{td}%
(\Sigma) = -\mathrm{rk}(\mathcal{E})\int_{\Sigma}c_{1}(\Sigma) = 2(g-1)\mathrm{rk}(\mathcal{E})
~~.
\]

It is also important to observe that
\begin{equation}
\label{h1split}H^{1}(\mathbf{M}_{(\eta,\sigma)}) \simeq\mathbb{H}%
^{1}(\mathbf{E}(\eta,\sigma)) \simeq\mathrm{coker}~H^{0}(\theta_{(\eta
,\sigma)}) \oplus\mathrm{ker}~H^{1}(\theta_{(\eta,\sigma)}) ~~,
\end{equation}
where $H^{0}(\theta_{(\eta,\sigma)}) ~:~ H^{0}(\mathcal{E}(\eta)) \to
H^{0}(\mathcal{E}(\eta)\otimes K_{\Sigma})$ and $H^{1}(\theta_{(\eta,\sigma)})
~:~ H^{1}(\mathcal{E}(\eta)) \to H^{1}(\mathcal{E}(\eta)\otimes K_{\Sigma})$
are the maps induced by $\theta_{(\eta,\sigma)}$ in cohomology; see the proof of
\cite[Proposition 3.1.11]{Bo} and also \cite{J1}.

The dual operator:
$$ \mathcal{D}_{(\eta,\sigma)}^{*} ~:~
L_{p}^{2}\left(\mathcal{E}(\eta)\otimes(\Lambda_{\Sigma}^{1,0}\oplus\Lambda_{\Sigma}^{0,1})\right)
\longrightarrow
L_{p-1}^{2}\left(\mathcal{E}(\eta)\otimes(\Lambda_{\Sigma}^{0}\oplus\Lambda_{\Sigma}^{1,1})\right) $$
is given by
$$ \mathcal{D}_{(\eta,\sigma)}^{*}(\varphi_1,\varphi_2) = 
(\theta^\vee_{(\eta,\sigma)}\varphi_1+\overline{\partial}^*_\eta\varphi_2 , -\overline{\partial}_\eta\varphi_1-\theta_{(\eta,\sigma)}\varphi_2) ~~.$$
The corresponding Laplacian operator
$$ \mathcal{D}_{(\eta,\sigma)}^{*}\mathcal{D}_{(\eta,\sigma)} ~:~
L_{p+1}^{2}\left(\mathcal{E}(\eta)\otimes(\Lambda_{\Sigma}^{0}\oplus\Lambda_{\Sigma}^{1,1})\right)
\longrightarrow
L_{p-1}^{2}\left(\mathcal{E}(\eta)\otimes(\Lambda_{\Sigma}^{0}\oplus\Lambda_{\Sigma}^{1,1})\right)
$$
is then given by
\begin{equation}\label{lapl}
\mathcal{D}_{(\eta,\sigma)}^{*}\mathcal{D}_{(\eta,\sigma)}(\psi_0,\psi_2) = 
(i\Lambda(\theta_{(\eta,\sigma)}^*\theta_{(\eta,\sigma)}+\partial_\eta\overline{\partial}_\eta)\psi_0 ,
i(\overline{\partial}_\eta\partial_\eta+\theta_{(\eta,\sigma)}\theta_{(\eta,\sigma)}^*)\Lambda\psi_2) ~~.
\end{equation}
To obtain this expression, we have used the K\"ahler identities
$\overline{\partial}_\eta^*=i[\partial_\eta,\Lambda]$ and 
$\partial_\eta^*=-i[\overline{\partial}_\eta,\Lambda]$, and $\overline{\partial}_\eta\theta_{(\eta,\sigma)}=0$.

Let us now denote by $G_{(\eta,\sigma)}$ the Green operator associated to the Laplacian $\mathcal{D}_{(\eta,\sigma)}^{*}\mathcal{D}_{(\eta,\sigma)}$ (i.e., a right inverse for $\mathcal{D}_{(\eta,\sigma)}^{*}\mathcal{D}_{(\eta ,\sigma)}$). We are then able to define
a family of operators
$$ J^{\vee}\times H^{0}(K_{\Sigma}) \ni(\eta,\sigma) \longmapsto P_{(\eta,\sigma)}:
L^{2}_{p}\left(\mathcal{E}\otimes(\Lambda_{\Sigma}^{1,0}\oplus
\Lambda_{\Sigma}^{0,1})\right)  \to L^{2}_{p}\left(\mathcal{E}\otimes(\Lambda_{\Sigma
}^{1,0}\oplus\Lambda_{\Sigma}^{0,1})\right) $$
by setting
\[
P_{(\eta,\sigma)}=1-\mathcal{D}_{(\eta,\sigma)}G_{(\eta,\sigma)}%
\mathcal{D}_{(\eta,\sigma)}^{*} ~~.
\]
It is straightforward routine to check that $P$ is a projector for
$\widehat{\mathcal{E}}:=\mathrm{coker}~\mathcal{D}\hookrightarrow H^{-}$. In other words,
that each $P_{(\eta,\sigma)}$ is idempotent, selfadjoint, restricts to the
identity on $\mathrm{coker}~\mathcal{D}_{(\eta,\sigma)}$ and vanishes on its
orthogonal complement.

Hence, having fixed the trivial connection
$\underline{d}$ on $H^{-}\to J^{\vee}\times H^{0}(K_{\Sigma})$, we may define a
unitary connection $\widehat{\nabla}$ on the bundle $\widehat{\mathcal{E}}\to J^{\vee
}\times H^{0}(K_{\Sigma})$ by demanding commutativity of the following
diagram:
\[
\xymatrix{
\Lambda^{0}_{J^{\vee}\times H^{0}(K_{\Sigma})}(\widehat{\mathcal{E}}) \ar[d] \ar[r]^{\widehat{\nabla}} & \Lambda^{1}_{J^{\vee}\times H^{0}(K_{\Sigma})}(\widehat{\mathcal{E}}) \ar[d] \\
\Lambda^{0}_{J^{\vee}\times H^{0}(K_{\Sigma})}(H^{-}) \ar[r]^{\underline{d}} & \Lambda^{1}_{J^{\vee}\times H^{0}(K_{\Sigma})}(H^{-})
}
\]

\begin{definition}
The vector bundle $\widehat{\mathcal{E}}\to J^{\vee}\times H^{0}(K_{\Sigma})$, endowed
with the connection $\widehat{\nabla}$, is called the \emph{Nahm transform} of $(\mathcal{\mathcal{E}},\theta)$.
\end{definition}

\begin{remark}\rm
We should also point out that the space $J^{\vee}\times H^{0}(K_{\Sigma})$ can be regarded as the moduli
space of rank one, degree zero Higgs bundles on $\Sigma$:
\[
J^{\vee}\times H^{0}(K_{\Sigma}) \ni(\eta,\sigma)\longleftrightarrow(L_{\eta},\sigma) ~~.
\]
Therefore, the Nahm transform presented here is a \emph{flat Nahm transform} in the sense of \cite{J5}.
\end{remark}

\begin{proposition}
The Nahm transforms of isomorphic, stable, degree zero Higgs bundles are gauge-equivalent.
\end{proposition}

\begin{proof}
Indeed, if $\mathbf{E}=\left(\ce,\theta\right)$ and $\mathbf{E}^{\prime}=\left(\ce^{\prime},\theta^{\prime}\right)$
are isomorphic as Higgs bundles, then there is a unitary bundle isomorphism $u:\mathbf{E}\to\mathbf{E}^{\prime}$
such that $\theta'=u\theta u^{-1}$. It follows that  $\left(\theta^{\prime}\right)^{\ast} = u \left(  \theta^{\ast}\right) u^{-1}$ since $u^{\ast}=u^{-1}$ and that
$$ \left(\overline{\partial^{\prime}}\right)^{\ast}\psi =
u \overline{\partial}^{\ast}\left(u^{-1}\psi\right) ~~ .$$
Following our recipe, we form $\mathcal{D}^{\prime}$, and notice that
$$ \left(u\otimes1\right) \mathcal{D}_{(\eta,\sigma)}^{\prime*}\psi =
\mathcal{D}_{(\eta,\sigma)}^{*} \left(\left(u\otimes1\right)^{-1}\psi\right) $$
i.e., $U=\left(u\otimes1\right)^{-1}$ defines a natural isomorphism
${\rm coker}~\mathcal{D}_{(\eta,\sigma)}^{\prime}\to {\rm coker}~\mathcal{D}_{(\eta,\sigma)}$
which is independent of $(\eta,\sigma)$ in the sense that $\underline{d}U=0$. Now the identities
\begin{align*}
\mathcal{D}_{(\eta,\sigma)}^{*}\mathcal{D}_{(\eta,\sigma)}G_{(\eta,\sigma)}\psi & = \psi \\
\mathcal{D}_{(\eta,\sigma)}^{\prime*}\mathcal{D}_{(\eta,\sigma)}^{\prime}G_{(\eta,\sigma)}^{\prime}\psi & =\psi
\end{align*}
imply $G_{(\eta,\sigma)}^{\prime}\psi = U G_{(\eta,\sigma)}\left(U^{-1}\psi\right)$, hence
$P_{(\eta,\sigma)}^{\prime}  = U P_{(\eta,\sigma)} U^{-1}$. Therefore we have
$$ \widehat{\nabla^{\prime}} = \left(UPU^{-1}\right)\underline{d} =
U P\underline{d} U^{-1} = \widehat{\nabla}^{U} ~~. $$
as we wished to show.
\end{proof}

\begin{proposition}
The Nahm transform is additive, i.e.:
\[
\left(  \widehat{\mathcal{E}_{1}\oplus \mathcal{E}_{2}},\widehat{\nabla_{1}\oplus\nabla_{2}%
}\right)  \simeq\left(  \widehat{\mathcal{E}_{1}}\oplus\widehat{\mathcal{E}_{2}},\widehat
{\nabla_{1}}\oplus\widehat{\nabla_{2}}\right)  ~~.
\]

\end{proposition}

It follows from this proposition that the Nahm transform is well defined for
solutions of Hitchin's equations that are \emph{without flat factors}, or
equivalently for polystable Higgs bundles of zero slope and no rank one summands.

\begin{proof}
In fact, if $\mathcal{E}$ is the direct sum of the bundles $\mathcal{E}_{1},\mathcal{E}_{2}$, it is fairly
obvious that, under the decomposition:
$$ L^2_p\left(\mathcal{E}(\eta)\otimes(\Lambda_{\Sigma}^{1,0}\oplus\Lambda_{\Sigma}^{0,1})\right) \simeq
L^2_p\left(\mathcal{E}_1(\eta)\otimes(\Lambda_{\Sigma}^{1,0}\oplus\Lambda_{\Sigma}^{0,1})\right) \oplus
L^2_p\left(\mathcal{E}_2(\eta)\otimes(\Lambda_{\Sigma}^{1,0}\oplus\Lambda_{\Sigma}^{0,1})\right) ~~,$$
the operator
$$ \mathcal{D}_{(\eta,\sigma)} :
L^2_p\left(\mathcal{E}(\eta)\otimes(\Lambda_{\Sigma}^{0}\oplus\Lambda_{\Sigma}^{1,1})\right) \to
L^2_{p-1}\left(\mathcal{E}(\eta)\otimes(\Lambda_{\Sigma}^{1,0}\oplus\Lambda_{\Sigma}^{0,1})\right) $$
splits as $\mathcal{D}_{1,(\eta,\sigma)}\oplus\mathcal{D}_{2,(\eta,\sigma)}$, hence
${\rm coker}~\mathcal{D}_{(\eta,\sigma)}={\rm coker}~\mathcal{D}_{1,(\eta,\sigma)}\oplus
{\rm coker}~\mathcal{D}_{2,(\eta,\sigma)}$. Moreover,
$G_{(\eta,\sigma)}=G_{1,(\eta,\sigma)}\oplus G_{2,(\eta,\sigma)}$ and
$P_{(\eta,\sigma)}=P_{1,(\eta,\sigma)}\oplus P_{2,(\eta,\sigma)}$
so that
$$ \widehat{\nabla} = P_{(\eta,\sigma)}\underline{d} = \widehat{\nabla_{1}} \oplus \widehat{\nabla_{2}} $$
as claimed.
\end{proof}

\begin{proposition} \label{pb}
The Nahm transform of an irreducible solution of Hitchin's equations with a vanishing Higgs field is the pullback from a bundle with unitary connection $(\widehat{\mathcal{E}},\widehat{\nabla})$ over $J^{\vee}$. Moreover, the fiber $\widehat{\mathcal{E}}_{\eta}$ can be identified with 
$H^{1}(\mathcal{E}^{\vee}(-\eta))^{\vee}\oplus H^{1}(\mathcal{E}(\eta))$ (where $^{\vee}$ here means taking duals).
\end{proposition}

\begin{proof}
Let $\left(\ce,\theta\right)$ be the corresponding stable Higgs bundle with slope zero and $\theta=0$.
Then the fibers of the transformed bundle $\widehat{\mathcal{E}}\to J^{\vee}\times H^{0}(K_{\Sigma})$
are given by the kernel of the operator
$\mathcal{D}_{(\eta,\sigma)}^*=-\overline{\partial}_{\eta}^{\ast}-\overline{\partial}_{\eta}+1\otimes\varphi_{\sigma}$,
where $\varphi_{\sigma}=\overline{\sigma}-\sigma$.
Now given $t$ in the interval $I=[0,1]$, consider the operator
$$ T_{(\eta,\sigma,t)} ~:~ L_{p}^{2}\left(\mathcal{E}\otimes(\Lambda_{\Sigma}^{0}\oplus\Lambda_{\Sigma}^{1,1})\right) \to
L_{p-1}^{2}\left(\mathcal{E}\otimes(\Lambda_{\Sigma}^{1,0}\oplus\Lambda_{\Sigma}^{0,1})\right) $$
$$ T_{(\eta,\sigma,t)} = \mathcal{D}_{(\eta,t\sigma)} ~~.$$
By our previous arguments, $T_{(\eta,\sigma,t)}$ is an injective Fredholm operator
for each $(\eta,\sigma,t)\in J^{\vee}\times H^{0}(K_{\Sigma})\times I$, with the same
index as $\mathcal{D}_{(\eta,\sigma)}$. So we have a complex vector bundle
$$ \widehat{\mathcal{E}_{I}}\longrightarrow J^{\vee}\times H^{0}(K_{\Sigma})\times I $$
with fibers given by ${\rm coker}~T_{(\eta,\sigma,t)}$, which gives us an isomorphism
of vector bundles between $\widehat{\mathcal{E}_{0}}=\widehat{\mathcal{E}_{I}}|_{J^{\vee}\times H^{0}(K_{\Sigma})\times\{0\}}$
and $\widehat{\mathcal{E}_{1}}=\widehat{\mathcal{E}_{I}}|_{J^{\vee}\times H^{0}(K_{\Sigma})\times\{1\}}$
But $\widehat{\mathcal{E}_{1}}$ coincides with the transformed bundle $\widehat{\mathcal{E}}$, while
the fibers of $\widehat{\mathcal{E}_{0}}$ are given by
$$ \left(\widehat{\mathcal{E}_{0}}\right)_{(\eta,\sigma)} = {\rm coker}~ \mathcal{D}_{(\eta,0)} ~~. $$
Therefore, $\widehat{\mathcal{E}_{0}}$ coincides with the pullback of the bundle
\begin{equation}\label{pbb}
J^{\vee}\ni\eta \longmapsto {\rm coker}~\left(\overline{\partial}_{\eta}+\overline{\partial}_{\eta}^{\ast}\right)
\end{equation}
by the obvious projection
$$ J^{\vee}\times H^{0}(K_{\Sigma})  \longrightarrow J^{\vee} ~~.$$
A similar argument allows us to conclude that the transformed connection $\widehat{\nabla}$ is also the
pullback from a connection on the bundle (\ref{pbb}). The final claim follows from (\ref{h1split}) and
Serre duality.
\end{proof}

In other words, the Nahm transform of an irreducible solution of Hitchin's equations with a vanishing
Higgs field coincides with the pullback of a (possibly trivial) extension of the Nahm transform
(in the sense of \cite[Section 4]{TP}) of $(\mathcal{E},\nabla)$ by the Nahm transform of  $(\mathcal{E}^{\vee},\nabla^{\vee})$.

Everything said so far is independent of the choice of a complex structure on $J^{\vee}\times H^{0}(K_{\Sigma})$.
We now fix a hyperk\"ahler structure $(I_{1},I_{2},I_{3})$ on $J^{\vee}\times H^{0}(K_{\Sigma})$, as defined in the end of Section \ref{s}.

\begin{proposition}\label{qi}
The transformed connection is hyperholomorphic.
\end{proposition}

In particular, if $g$ is even (so that $J^{\vee}$ is hyperk\"ahler) and the Higgs field vanishes, then the transformed connection is the pullback of a hyperholomorphic connection on $J^\vee$.

\begin{proof}
One can think of the complex (\ref{c1}) as a monad (of trivial, infinite dimensional vector bundles)
over the hyperk\"ahler manifold $J^{\vee}\times H^{0}(K_{\Sigma})$, so that our transformed bundle is
exactly the cohomology of this monad. Since the operators in (\ref{c1}) vary holomorphically with respect
to any of the complex structures induced by the hyperk\"ahler structure, we conclude from general theory
(cf. \cite[Section 3.1]{DK}) that $\widehat{\mathcal{E}}$ has a holomorphic structure with respect to each complex
structure in $J^{\vee}\times H^{0}(K_{\Sigma})$, with which $\widehat{\nabla}$ is compatible. Hence
$F_{\widehat{\nabla}}$ must be of type $(1,1)$ with respect to all complex structures.
\end{proof}

Summing up the work done so far, we have proved the first part of our main theorem, namely that the Nahm
transform of a solution of Hitchin's equations on a Riemann surface $\Sigma$ is a Hermitian vector bundle
over $J^{\vee}\times H^{0}(K_{\Sigma})$ provided with a unitary connection which is hyperholomorphic.
The second part of the Main Theorem is proved below.

\begin{proposition}
The holomorphic structure induced by $\widehat{\nabla}$ on $\widehat{\mathcal{E}}$ with
respect to a product complex structure extends to a holomorphic bundle over
$J^{\vee}\times\mathbb{P}(H^{0}(K_{\Sigma})\oplus\mathbb{C})$.
\end{proposition}

\begin{proof}
Consider a product complex structure on $J^\vee\times H^0(K_\Sigma)$, and let $\widehat{\ce}$ denote
the holomorphic vector bundle over $J^\vee\times H^0(K_\Sigma)$ given by the transformed
bundle $\widehat{E}$ provided by the holomorphic structure
$\overline{\partial}_{\widehat{\nabla}}$. We must exhibit a holomorphic vector bundle
$\overline{\hat\ce}$ over $J^\vee\times\pp^g$ such that
$\overline{\hat\ce}|_{J^\vee\times\C^g} \simeq \widehat{\ce}$,
where we have identified $\pp(H^0(K_\Sigma)\oplus\C)$ with $\pp^g$ and $H^0(K_\Sigma)$
with $\C^g$, the big open cell in $\pp^g$.
Following Bonsdorff \cite{Bo}, we fix a basis of sections $z_{0},...,z_{g}$ of
$H^{0}(\mathcal{O}_{\pp^g}(1))$ or, in more pedestrian terms, homogeneous coordinates
for $\pp^g$, where $H^0(K_\Sigma)$ corresponds to the affine space $z_{0}\neq0$.
Let then $\left\{\sigma_{1},...,\sigma_{g}\right\}$ be a basis of holomorphic $1$-forms
on $\Sigma$, and define the Higgs field
$$ \theta_{(\eta,[z_{0}:...:z_{g}])} =
z_{0}\left( \theta\otimes\id_{L_{\eta}} \right) + \id_{E(\eta)}\otimes\sigma ~~, $$
so that
$$ \theta_{(\eta,[1:z_{1}:...:z_{g}])} = \theta_{(\eta,\sigma)} ~~{\rm and}~~
\theta_{(\eta,[0:z_{1}:...:z_{g}])} = \id_{E(\eta)}\otimes\sigma ~~, $$
where $\sigma=\sum_{k=1}^{g}z_{k}\sigma_{k}$.
Thus we obtain a holomorphic family of Higgs bundles
$$ \mathbf{E}_{(\eta,[z])} =
\left\{ \ce(\eta) \overset{\theta_{(\eta,[z])}}{\longrightarrow} \ce(\eta)\otimes K_\Sigma \right\} $$
parameterized by $J^\vee\times\pp^g$. For $(\eta,[z])$ on the affine subset $J^{\vee}\times H^{0}(K_{\Sigma})$
it follows from Proposition \ref{v1} that $\mathbb{H}^{0}\left(\mathbf{E}_{(\eta,[z])}\right)=
\mathbb{H}^{2}\left(\mathbf{E}_{(\eta,[z])}\right)=0$. Thus we obtain a holomorphic vector bundle
$$ J_{\Sigma}^{\vee}\times H^{0}(K_{\Sigma})\ni(\eta,[0:z_{1}:...:z_{g}])\longmapsto
\mathbb{H}^{1}\left(\mathbf{E}_{(\eta,[0:z_{1}:...:z_{g}])}\right)
$$ which, by general theory (cf. \cite[Section 3.1]{DK}), coincides
with $\widehat{\ce}$. Now we claim that the hypercohomology
$\mathbb{H}^{\bullet}\left(\mathbf{E}_{(\eta,[0:z_1:\dots:z_g])}\right)$
is also concentrated in dimension $1$. It will then follow that the
assignment
$$ J^{\vee}\times\pp^g\ni(\eta,[z]) \longmapsto \mathbb{H}^{1}\left(\Sigma,\mathbf{E}_{(\eta,[z])}\right) $$
provides a holomorphic extension $\overline{\hat\ce}$ of $\widehat{\ce}$ to $J^{\vee}\times\pp^g$.

Indeed, $\mathbb{H}^{0}\left(\mathbf{E}_{(\eta,[0:z_1:\dots:z_g])}\right)$ can
be identified with global holomorphic sections of the kernel of the bundle morphism
$$ \ce(\eta) \overset{\id_{\ce(\eta)}\otimes\sigma}{\longrightarrow} \ce(\eta)\otimes K_{\Sigma} ~~. $$
Since $\sigma=\sum_{k=1}^{g}z_{k}\sigma_{k}$ is nontrivial and vanishes only at finitely many points, we
conclude that $\ker\left(1\otimes\sigma\right)$ is the zero sheaf, hence
$\mathbb{H}^{0}\left(\mathbf{E}_{(\eta,[0:z_1:\dots:z_g])}\right)=0$.
Now Serre duality allows us to identify
\[
\mathbb{H}^{2}\left(\mathbf{E}_{(\eta,[0:z_1:\dots:z_g])}\right) \simeq\mathbb{H}^{0}\left(\mathbf{E}_{(\eta,[0:z_1:\dots:z_g])}^{\vee}\otimes K_{\Sigma}\right)^{\vee}
\]
which again we regard as the space of global holomorphic sections of the kernel of the morphism
$$ \ce(\eta)^{\vee} \overset{\id_{\ce(\eta)^{\vee}}\otimes\sigma}{\longrightarrow}
\ce(\eta)^{\vee}\otimes K_{\Sigma} ~~. $$
Arguing as above, we conclude that $\mathbb{H}^{2}\left(\mathbf{E}_{(\eta,[0:z_1:\dots:z_g])}\right)=0$.
\end{proof}

\begin{remark}\rm
In principle, the Nahm transformed pair $(\widehat{\ce},\widehat{\nabla})$ depends on the
choice of point $p\in\Sigma$, which fixes the Abel map. It would be interesting to analyze exactly how
this dependence goes. One idea would be to define a ``universal" transformed bundle
$\widehat{\mathbb{E}}\to\Sigma\times J^\vee\times H^0(K_\Sigma)$, such that
$\widehat{\mathbb{E}}|_{\{p\}\times J^\vee\times H^0(K_\Sigma)}$ coincides with the Nahm transform
of $(\ce,\theta)$ induced by the point $p$.
\end{remark}

\end{document}